\documentclass[a4paper,11pt]{article}

\usepackage{graphicx,epstopdf,epsfig,color}
\usepackage{amsmath,amssymb}

\newcommand{\dd}{{\mathrm d}}

\newcommand{\eE}{{\mathrm e}}

\newcommand{\cI}{\mathcal{I}}

\newcommand{\un}[1]{\ \textrm{#1}}

\newcommand{\be}{\begin{equation} }
\newcommand{\ee}{\end{equation}}
\newcommand{\bmath}{\begin{displaymath} }
\newcommand{\emath}{\end{displaymath} }
\newcommand{\lab}[1]{\label{#1}}

\newcommand{\half}{\mbox{$\frac{1}{2}$}}

\newcommand{\I}[1]{{\, \dd {#1}}}

\newcommand{\K}{{\,\mathrm K}}

\newcommand{\ijh}[1]{{#1}}
\newcommand{\ijhdel}[1]{}

\title{A Mathematical Model for Flash Sintering}

\author{I.J. Hewitt \\
Mathematical Institute, University of Oxford, \\
Woodstock Road, Oxford, OX2 6GG, UK,\\
\\
A.A. Lacey \\
Maxwell Institute for Mathematical Sciences and
Department of Mathematics,\\
School of Mathematical and Computer Sciences, Heriot-Watt University, \\
Riccarton, Edinburgh, EH14 4AS, UK,\\
\\
and \\
\\
R.I. Todd \\
Department of Materials, University of Oxford, \\
Parks Road, Oxford, OX1 3PH, UK}
\date{}
\begin{document}

\maketitle

\begin{abstract}

A mathematical model is presented for the Joule heating that occurs
in a ceramic powder compact during the process of flash sintering.
The ceramic is assumed to have an electrical conductivity that increases
with temperature, and this leads to the possibility of runaway heating
that could facilitate and explain the rapid sintering seen in experiments.
We consider reduced models that are sufficiently simple to enable
concrete conclusions to be drawn about the mathematical nature
of their solutions. In particular we discuss how different local
and non-local reaction terms, which arise from specified
experimental conditions of fixed voltage and current,
lead to thermal runaway or to stable conditions.
We identify incipient thermal runaway as a necessary condition
for the flash event, and hence identify the conditions under
which this is likely to occur.

\end{abstract}


%*******************************************************************

%%%%%%%%%%%%%%%%%%%%%%%%%%%%%%%%%%%%%%%%%%%%%%%%%%%%%%%%%%%%%%%%

\section{Introduction}           \label{sec:intro}

%%%%%%%%%%%%%%%%%%%%%%%%%%%%%%%%%%%%%%%%%%%%%%%%%%%%%%%%%%%%%%%%

Flash sintering is a novel method for sintering ceramic materials,
\ijh{performed by} simultaneously heating a powder compact in a furnace
while passing through it an electric current \cite{Coletal10}.
At a critical furnace temperature or applied voltage,
there is an electrical power spike accompanied by rapid sintering
in a matter of seconds.  Sintering usually occurs through
grain-boundary diffusion, which transports matter from grain boundaries
into pore spaces and leads to densification of the material.
Its rate depends strongly on temperature, and conventional sintering
({\it i.e.}\ without the electric current) can take several hours
at furnace temperatures $\sim 1450^{\circ}$C \cite{Coletal10}.
By contrast, flash-sintering can occur in several seconds at
furnace temperatures $\sim 850^{\circ}$C. It has been estimated
\cite{Raj12} that a temperature of $\sim 1900^{\circ}$C would be
required to enable such rapid conventional sintering. 

There is still uncertainty over how exactly the imposition
of the electric field enables such rapid sintering.
Joule heating can raise the sample temperature above the furnace temperature,
and if sufficiently high temperatures were reached, the process might
be explained simply by the enhanced sintering rate at higher temperatures.
However, Joule heating is counteracted by radiative cooling,
and some authors have suggested that it would be ineffective at
increasing the temperature to the required levels \cite{Coletal10,Raj12}.
\ijh{It has therefore been suggested that local heating at the grain boundaries may be responsible
\cite{Fraetal12,FraRaj13}, or that the electric field leads
to new mass transport mechanisms such as greatly
increased concentrations of vacancies or interstitials,
\cite{NSR}.}

On the other hand, there is considerable difficulty in measuring or
estimating the sample's actual temperature, which can vary rapidly
during the process.  It is possible that Joule heating could have
a much more considerable effect than has been considered up to now
if one accounts for the reduction
in resistivity that is seen to accompany sintering.

The purpose of this paper is to provide a mathematical model
for the Joule heating of sintering material
that can help us to understand the flash sintering process.
In order to provide some concrete statements about this factor alone,
we adopt a number of simplifications; in particular we ignore
any material changes that occur as a result of the sintering,
and concentrate solely on the effect of an assumed temperature-dependence
of the electrical conductivity.  Specifically, we assume
an Arrhenius form for the conductivity.  The result is a mathematical model
for the electric field and the temperature in the sample,
with many similarities to the so-called thermistor problem
\cite{Cim90,Fowetal92}. 

The model is formulated in \S \ref{sec:model} Then
we study a number of reductions, in which the thermal problem is expressed
as a semi-linear reaction diffusion problem with either a local
or a non-local reaction term in \S \ref{sec:reduced} Two of the
reductions are obtained in limiting cases of negligible heat loss, from either
the ends (electrodes) or the sides. This allows us to cover extreme
cases, each of which might be relevant in practice
because of different experimental
set-ups or because of poorly constrained physical parameters. We make use of a number
of existing results concerning the solutions to these problems
to find when rapid temperature rise will occur, and to describe
the nature of the temperature evolution.  We show that under
a sufficiently high specified voltage, both local
and non-local reaction terms can lead to `blow-up' behaviour,
but that under a specified current  the temperature tends towards
a globally-attracting steady state.  Since flash sintering experiments
generally involve a switch from voltage control to current control
(to prevent runaway power consumption), a rapid rise in temperature
is followed by relaxation towards equilibrium, as has been observed
in numerical models of the process \cite{Graetal11,Todd}.
We close with a discussion of how these models
might be extended to help further understand the flash-sintering process in \S \ref{sec:concl}
Quantitative estimates suggest that Joule heating may be sufficient
to explain the flash sintering process, {at} least in some circumstances.

%%%%%%%%%%%%%%%%%%%%%%%%%%%%%%%%%%%%%%%%%%%%%%%%%%%%%%%%%%%%%%%%

\vspace*{0.5cm}
\setcounter{equation}{0}
\section{Modelling the Process}          \label{sec:model}

%%%%%%%%%%%%%%%%%%%%%%%%%%%%%%%%%%%%%%%%%%%%%%%%%%%%%%%%%%%%%%%%

\subsection{The dimensional model}      \label{subsec:dimensional}

% % % % % % % % % % % % % % % % % % % % % % % % % % % % % % % % % % % % 

\begin{figure}
\centering
\includegraphics[width=0.5\textwidth]{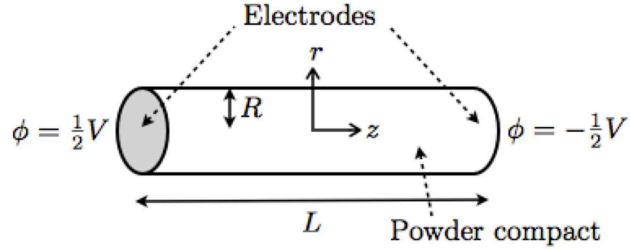}

\caption{\label{schematic}Model set-up considered in this paper.  The powder compact is held between two electrodes, at electric potentials $\pm \half V$.  Radiative and convective heat loss occurs from the free sides of the cylindrical compact, and from the electrodes. }
\end{figure}

We consider an axisymmetric sample occupying
the region $0< r < R(z)$, $-\half L<z< \half L$,
see Figure~\ref{schematic}.

%%%%%%%%%%%%%%%%%%%%%%%%%%%%%%%%%%%%%%%%%%%%%%%%%%

%%%%%%%%%%%%%%%%%%%%%%%%%%%%%%%%%%%%%%%%%%%%%%%%%%

The electrical conductivity of the ceramic $\sigma(T)$
is a function of temperature $T$, so the electric potential
$\phi$ satisfies
\be \label{phieq}
\nabla \cdot ( \sigma(T) \nabla \phi ) = 0,
\ee
with 
\be \label{phibc1}
\phi_r - R_z \phi_z = 0 \quad {\rm at} \quad r = R(z),
\ee
\be \label{phibc2}
\phi = \mp \half V \quad {\rm at} \quad z = \pm \half L. 
\ee
The current is 
\be \label{cur}
I = - 2\pi \int_{0}^{R} \sigma(T) \phi_z\, r \I r,
\ee
and is seen from (\ref{phieq}) to be independent of $z$. We prescribe
\emph{either} the voltage $V(t) = V_0$ \emph{or}
the current $I(t) = I_0$;
in the latter case the voltage is determined from
the integral constraint (\ref{cur}). The experimental procedure usually involves applying
a fixed voltage $V_0$ initially, but once a maximum current $I_0$
is reached, the current is maintained at this value
and the voltage is subsequently allowed to vary.
{Note that contact resistances, between the electrodes
and the sample, have been neglected. They are expected to be small
and the good agreement between the predictions of this model and
the experimental results of \cite{Todd} would suggest that neglecting
them, as has been done by other authors, is reasonable.
We should also note that this is the first
stage in our investigation and such resistances could be easily
included in more sophisticated models. However, their values are
not known and to be included accurately these would have to be
measured.}

The thermal problem is
\be \label{Teq}
\rho c T_t = \nabla \cdot ( k \nabla T ) + \sigma(T) |\nabla \phi|^2,
\ee
with
\be \label{Tbc1}
k\frac{T_r - R_z T_z}{\sqrt{1+R_z^2}} = - \ijh{h_s}(T-T_{\infty}) - \epsilon S \left[T^4-T_{\infty}^4\right] \quad {\rm at} \quad r = R(z),
\ee
\be \label{Tbc2}
\pm T_z = -h_{e} (T-T_{\infty}) \quad {\rm at} \quad z = \pm \half L \, ,
\quad \mbox{ and } \quad T = T_\infty \, \mbox{ at } \, t = 0.
\ee
Here $\rho$ is the density of the ceramic powder,
$c$ is the specific heat capacity,
and $k$ is the thermal conductivity. These are treated as constants.
Cooling from the sides
of the sample is taken to be through a combination of radiation,
with emissivity $\epsilon$, \ijh{Stefan-Boltzmann} constant $S$,
and furnace temperature $T_{\infty}$, and conduction/convection,
which is parameterized by a heat transfer coefficient \ijh{$h_s$}.
Cooling from the ends (at the electrodes) is described by
a heat transfer coefficient $h_e$. The last condition in (\ref{Tbc2}) assumes
that the sample is initially at the furnace temperature
everywhere.

The electrical conductivity is an increasing function of temperature.
We assume the Arrhenius law,
\be \label{arr1}
\sigma(T) = A \eE^{-E/R_{g}T},
\ee 
where $E$ is the activation energy and $R_g$ the universal gas constant.

% % % % % % % % % % % % % % % % % % % % % % % % % % % % % % % % % % % % 

\subsection{The dimensionless model}      \label{subsec:nondim}

% % % % % % % % % % % % % % % % % % % % % % % % % % % % % % % % % % % % 

{Our axisymmetric models cover \ijh{some of} the
bone-shaped samples used in many experiments (see, for example,
\cite{Coletal10}). However, because of the preliminary nature
of the present investigation, here
we concentrate on the case of cylindrical symmetry
so that there is uniform radius $r = R$, as in the experiments of
\cite{GSPCJLR} and \cite{Todd2}. It is intended
that situations with $R$ varying with $z$ will be considered in
a forthcoming work.}

The variables are non-dimensionalised by putting
\be
z = L\hat{z}, \quad
r = R \hat{z}, \quad
t = t_0 \hat{t}, \quad
T = T_{0} + \Delta T \;\! {\theta}, \quad
\nonumber
\ee
\be
\phi = V_0 \hat{\phi}, \quad
V = V_0 \hat{V}, \quad
I = \frac{\sigma_0 V_0 \pi R^2 }{L} \hat{I}, \quad
\sigma(T) = \sigma_0 \hat{\sigma}(\theta), \quad
\ee
where we choose the thermal-conduction time scale $t_0 = {\rho c R^2}/{k}$,
the electrical-conductivity scale $\sigma_0 = A \eE^{-E/R_gT_0}$,
a sensible choice for $T_0$ is the furnace temperature $T_{\infty}$,
and we take the temperature perturbation scale
$\Delta T = R_g T_0^2 / E$ so that the dimensionless conductivity is
\be \label{nd2sig}
\hat{\sigma}(\theta) = \eE^{\theta/(1+\nu \theta)}, \qquad \nu = \frac{R_g T_0}{E}.
\ee
The parameter $\nu$ is typically small, so we make the approximation
$\nu \to 0$, as is standard for treatments of chemical reactions
(the Frank-Kamenetskii approximation).

Dropping hats, the resulting model for the electric field is
\be \label{ndphieq}
\delta^2 \left[ \eE^{\theta}\phi_z\right]_z +\frac{1}{r} \left[r\eE^{\theta}\phi_r\right]_r = 0,
\ee
\be \label{ndphibc1}
\phi_r = 0 \quad {\rm at} \quad r = 1,
\ee
\be \label{ndphibc2}
\phi = \mp \half V \quad {\rm at} \quad z = \pm \half,
\ee
and for the temperature,
\be \label{ndTeq}
\theta_t = \delta^2 \theta_{zz} +\frac{1}{r}\left[ r\theta_r \right]_r  + \lambda \eE^{\theta} \left( \phi_z^2+\frac{1}{\delta^2}\phi_r^2 \right),
\ee
\be \label{ndTbc1}
\theta_r  = - \ijh{\beta} \theta \quad {\rm at} \quad r = 1,
\ee
\be \label{ndTbc2}
\theta_z = \mp \ijh{\alpha} \theta \quad {\rm at} \quad z = \pm \half \, ,
\, \mbox{ and } \quad \theta = 0 \, \mbox{ at } \, t = 0.
\ee
Note that as well as neglecting terms in $\nu$ from the exponential
in (\ref{ndTeq}), terms in $\nu$
have also been dropped from (\ref{ndTbc1}).

The dimensionless current and its maximum value are given by 
\be
I = - 2 \int_{0}^{R} \eE^{\theta} \phi_z\, r \I r \le \mathcal{I} ,
\qquad  \mathcal{I} = \frac{I_0 L}{ \sigma_0 V_0  \pi R^2 }.
\ee
The other dimensionless parameters in the model are:
\be
\delta = \frac{R}{L}, \qquad
\lambda = \frac{\sigma_0 V_0^2 R^2}{k \Delta T L^2}, \qquad
%\alpha = \frac{h R}{k}, \qquad
\ijh{\beta = \frac{h_sR + 4 \epsilon S T_{0}^3 R}{k},\qquad
\alpha = \frac{h_e L}{k}}
\ee
representing the aspect ratio, the strength of the ohmic heating,
the combined strength of conductive and radiative cooling from the sides,
and strength of cooling from the ends.  Typical values are shown in
Table~\ref{parametervalues}, although there is significant uncertainty
in the values for experimental samples, particularly in
the appropriate parameters for the Arrhenius law.
Note that the dominant contribution to cooling from the sides
comes from the radiation term.

\begin{table}[t] 
\centering
\begin{tabular}{l l l}
\hline
$\rho$ & Sample density & $ 6050 \un{kg}\un{m}^{-3}$\\
$c$ & Specific heat & $ 600 \un{J}\un{kg}^{-1}\un{K}^{-1}$\\ 
$k$ & Thermal conductivity & $ 2.7 \un{J}\un{m}^{-1}\un{s}^{-1}\un{K}^{-1}$\\
$\epsilon$ & Emissivity & $0.7$\\
$S$ & Stefan-Boltzmann constant & $5.67 \times 10^{-8}\un{W}\un{m}^{-2}\un{K}^{-4}$\\
\ijh{$h_s$} & Side heat transfer coefficient & $10 \un{J}\un{m}^{-2}\un{s}^{-1}\un{K}^{-1}$\\
$h_e$ & Electrode heat transfer coefficient &  \ijh{$10 \un{J}\un{m}^{-2}\un{s}^{-1}\un{K}^{-1}$}\\
$A$ & Arrhenius rate factor & $ 9.3\times 10^5 \un{S}\un{m}^{-1}$ \\
$E$ & Activation energy & $171 \un{kJ}\un{mol}^{-1}$\\
$R_g$ & Gas constant & $8.31\un{J}\un{K}^{-1}\un{mol}^{-1}$\\
\hline
\end{tabular} 
\vspace{0.2cm}

\begin{tabular}{l l l}
\hline
$L$ & Sample length & $ 10 \un{mm}$\\ 
$R$ & Sample radius & $ 1.5 \un{mm}$\\ 
$T_{\infty}$ & Furnace temperature &  {$1110\,\K$} \\
$V_0$ & Initial voltage & $300 \un{V}$ \\
$I_0$ & Current limit & $ 0.5 \un{A}$ \\
$\sigma_0$ & Conductivity at $T_{\infty}$ & $ 8.30 \times 10^{-3} \un{S}\un{m}^{-1}$ \\
$t_0$ & Time scale & $3.025\un{s}$ \\

\hline
\end{tabular}
\hspace{0.85cm}
\begin{tabular}{l l}
\hline
$\delta$ & $0.15$ \\
$\lambda$ & $0.104$ \\
%$\alpha$ & $5.56 \times 10^{-3}$ \\
$\beta$ & $0.126$ \\
\ijh{$ \alpha$} & \ijh{$0.037$} \\
$\mathcal{I}$ & \ijh{$284$} \\
$\nu $ & $0.054$ \\
& \\
\hline
\end{tabular}
\caption{\lab{parametervalues}Typical parameter values and scales used for non-dimensionalisation, along with the corresponding dimensionless parameters. Experimental values for 3YSZ (3~mol~\% yttria stabilised zirconia) are from \cite{Todd}, along with their Arrhenius fit for the conductivity.  An appropriate value for the electrode heat transfer coefficient may vary enormously; we take the same as for the sides for illustration.}
\end{table}

\

General mathematical results on existence and blow-up of solutions
of the coupled equations (\ref{phieq}) and (\ref{Teq}),
or equivalently (\ref{ndphieq}) and (\ref{ndTeq}), subject
to more limited boundary conditions, can be found in, for
example, \cite{AC}. However, our approach here will be
to simplify the model, given the sizes of parameters in
Table~\ref{parametervalues}, so that we can be more
specific about the qualitative behaviour of temperature
and electric current and obtain relatively
simple criteria for the occurrence of flash. In particular, we consider
two limiting cases in which the problem is effectively reduced
to one spatial dimension.  First, we take the case
of thermally insulating electrodes, $\alpha = 0$, in which case
$\theta = \theta(r,t)$. Secondly, we take the case of weakly cooled sides,
$\beta = 0$, in which case $\theta = \theta(z,t)$.
Next we consider the case (as for the estimated values in
Table~\ref{parametervalues}) when the heating $\lambda$,
as well as the cooling rates $\alpha$ and $\beta$, are all small.
Finally, we look at the high-aspect-ratio limit, $\delta$ small.

%%%%%%%%%%%%%%%%%%%%%%%%%%%%%%%%%%%%%%%%%%%%%%%%%%%%%%%%%%%%%%%%

\vspace*{0.5cm}
\setcounter{equation}{0}
\section{Reduced Models}              \label{sec:reduced}

%%%%%%%%%%%%%%%%%%%%%%%%%%%%%%%%%%%%%%%%%%%%%%%%%%%%%%%%%%%%%%%%

\subsection{Thermally insulating electrodes}   \label{subsec:insulate}

% % % % % % % % % % % % % % % % % % % % % % % % % % % % % % % % % % % % 

{If the electrodes % are embedded within
% the sample and thermal conduction along them can be ignored,
remove heat slowly, or heat transfer between them and the sample is poor,
we can regard them as providing thermal insulation to the ends.
This corresponds to taking $\alpha = 0$, and indeed the small value of $\alpha$
in Table \ref{parametervalues} suggests this as a possible simplification.}
In that case there is nothing to induce
any $z$-dependence of the temperature so we have
$\theta = \theta(r,t)$, $\phi = \phi(z,t)$, and the aspect ratio
$\delta$ is removed from the problem.
Integrating (\ref{ndphieq}) we obtain
\be
\phi = -V z, \qquad I = 2 V \int_{0}^{1} \eE^{\theta}\, r \I r.
\label{eq:VandI} \ee
The heating term becomes $V^2\eE^{\theta}$, and, remembering that
either $V=1$ or $I = \mathcal{I}$, we have
\be \label{Veq}
V = {\rm min} \left( 1\, ,{\mathcal{I}} \Big/ \left( 2 \int_{0}^{1} \eE^{\theta}
r \I r \right) \right) .
\ee

Thus the problem for the temperature is
\be \label{nd2Teq}
\theta_t = \frac{1}{r}\left[ r\theta_r \right]_r  + \lambda \eE^{\theta},
\ee
\be \label{nd2Tbc1}
\mbox{with } \, \theta_r  = - \beta \theta \quad {\rm at} \quad r = 1 \, ,
\, \mbox{ and } \quad \theta = 0 \, \mbox{ at } \, t = 0,
\ee
while $2 \int_{0}^{1} \eE^{\theta}
r \I r \le \mathcal{I}$. If and when the current limit is reached,
the equation switches to
\be \label{nd2Teq2}
\theta_t = \frac{1}{r}\left[ r\theta_r \right]_r  +
\frac{\lambda \cI^2 \eE^{\theta}}{4 (\int_{0}^{1} \eE^{\theta} r \I r)^2} \, ,
\ee
with the same boundary conditions.  
Apart from the switch between voltage and current control,
these are standard models; to start with, a semi-linear
reaction diffusion problem, and later, a similar problem
but with a non-local reaction term. Many results are available
concerning the existence of steady states and blow-up solutions.

\begin{figure}
\centering
\includegraphics[width=0.49\textwidth]{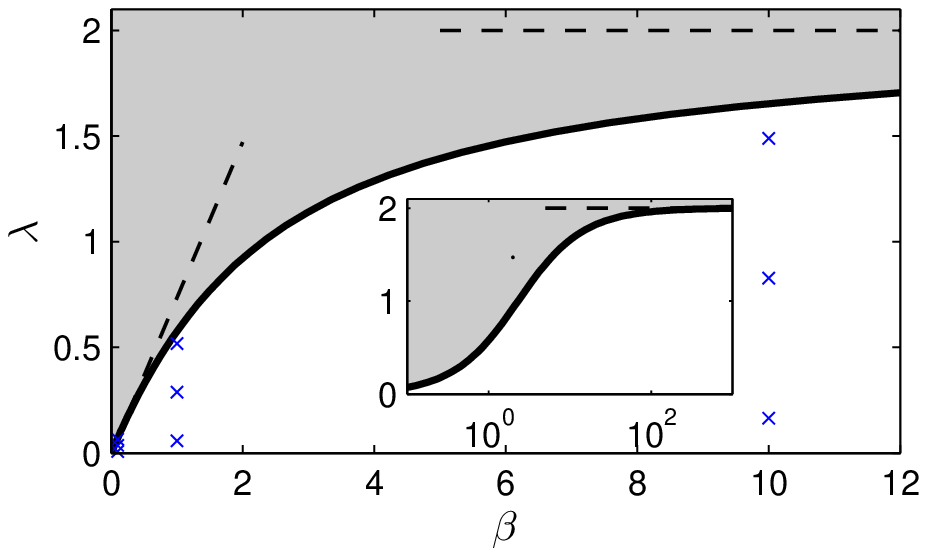}
\includegraphics[width=0.49\textwidth]{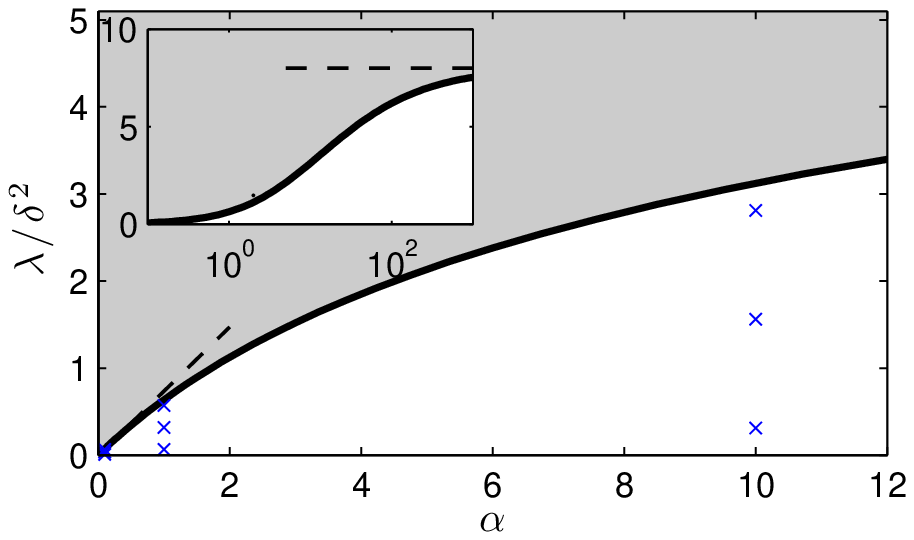}
\caption{\label{flashd_ststV2}\ijh{(a) Critical value $\lambda_c(\beta)$ for existence of a steady-state solution to (\ref{nd2Teq})-(\ref{nd2Tbc1}).  Dashed lines show the asymptotic behaviour and the inset shows the same on a logarithmic scale. (b) Critical value $\Lambda_c(\alpha)$ for existence of a steady-state solution to (\ref{redeq2})-(\ref{redeqbc}).} }
\end{figure}

In particular, for the fixed voltage problem {(\ref{nd2Teq})-(\ref{nd2Tbc1})},
there is a critical value $\lambda_c$ of $\lambda$ above which
there is no steady solution, {\cite{Amann, KC}}.
The critical value can be found from the exact solution,
\begin{equation}
\theta(r) = -2 \ln (c-b + br^2) \, ,
\label{eq:exactstst:1}
\end{equation}
where constants $b$ and $c$ are related by
\begin{equation}
{ \lambda = 8b(c-b) \, \, \mbox{ and } \,
b = \half \beta c(- \ln c).}
\label{eq:exactstst:2}
\end{equation}
Combining these conditions requires
\begin{equation} 
\lambda = 4 \beta c^2 (- \ln c)(1 + \half\beta \ln c) \, ,
\label{eq:exactstst:3}
\end{equation}
which has a solution for $c$ only if $\lambda < \lambda_c(\beta)$
(see also {\cite{Gelfand, JL}}). This curve is shown in
Figure~\ref{flashd_ststV2}.  In particular, $\lambda_c
\to 2$ for $\beta \to \infty$, while $\lambda_c \sim 2\beta/\eE$
for $\beta \to 0$. 

%%%%%%%%%%%%%%%%%%%%%%%%%%%%%%%%%%%%%%%%%%%%%%%%%%

%%%%%%%%%%%%%%%%%%%%%%%%%%%%%%%%%%%%%%%%%%%%%%%%%%

For values of $\lambda$ greater than $\lambda_c$, the solution
to {(\ref{nd2Teq})-(\ref{nd2Tbc1})} blows up
at a finite time, \cite{Lacey1}, and with this particular
initial data it does so at the centre $r=0$, \cite{FM},
and in a manner that leads to $\int\eE^\theta r\I r$ becoming
unbounded, {\cite{BB, HV1, HV2}}. This means,
from (\ref{eq:VandI}), that the current $I$ would become large,
and that the switch to the current-controlled problem
would necessarily occur before blow-up occurs.
{With large values of the limiting current, so that under current
control the sample is rapidly sintered, the incipient blow-up
behaviour can be identified as both
a necessary and a sufficient condition for the flash process
in the constant-voltage phase of a typical experiment,
since the rapid heating guarantees that the high temperatures
are reached and is also needed to generate them.}
(For more related results on problems such as
{(\ref{nd2Teq})-(\ref{nd2Tbc1})}, see \cite{BE}).
In \S~\ref{sec:concl} we discuss the interpretation of these conclusions,
and how they may be related to experimental parameters.

The problem under current control, (\ref{nd2Tbc1})-(\ref{nd2Teq2}),
is not so well studied, since it has a non-local term,
but it is known to have a solution
which is global in time, and tends to a unique steady state,
\cite{BL} (see also {\cite{Lacey3, Lacey4, Freitas}}).
The steady state is given by the same formula as above, (\ref{eq:exactstst:1}), but now with
\begin{equation}
\lambda \cI^2 = \frac {4 \beta (-\ln c)}{c^2(1 + \half\beta\ln c)} \quad
\mbox{ for } \, \eE^{-2/\beta} < c < 1  \, .
\label{eq:exactstst:4}
\end{equation}
Contrary to the voltage-controlled problem, this has a solution for $c$
regardless of the value of $\lambda$, so the steady state always exists.

\begin{figure}
\centering
\includegraphics[width=0.32\textwidth]{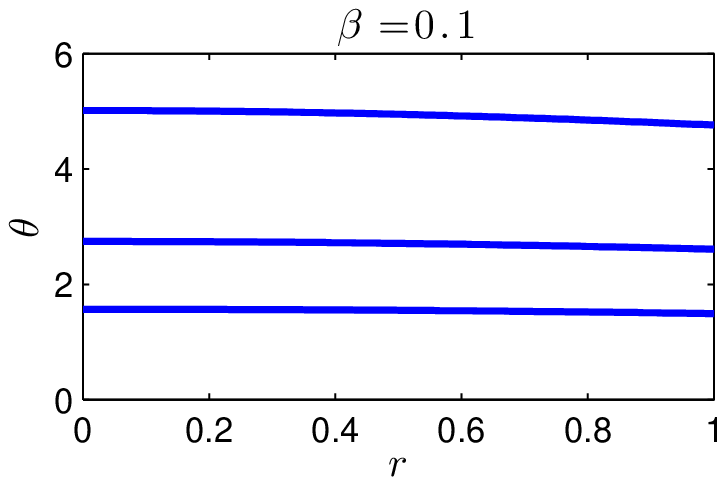}
\includegraphics[width=0.32\textwidth]{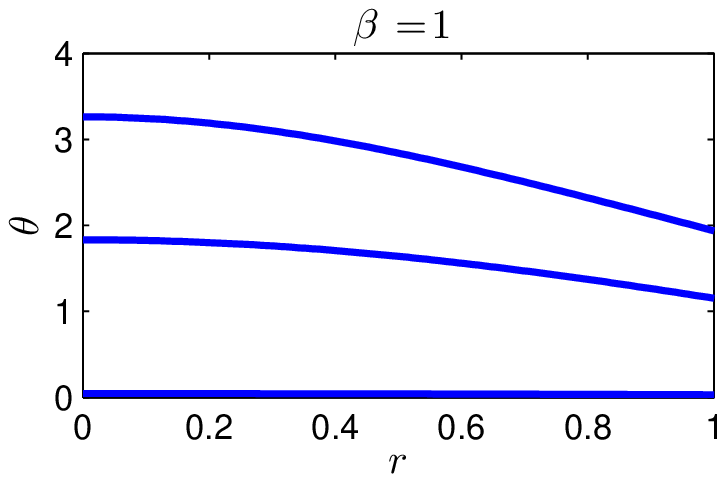}
\includegraphics[width=0.32\textwidth]{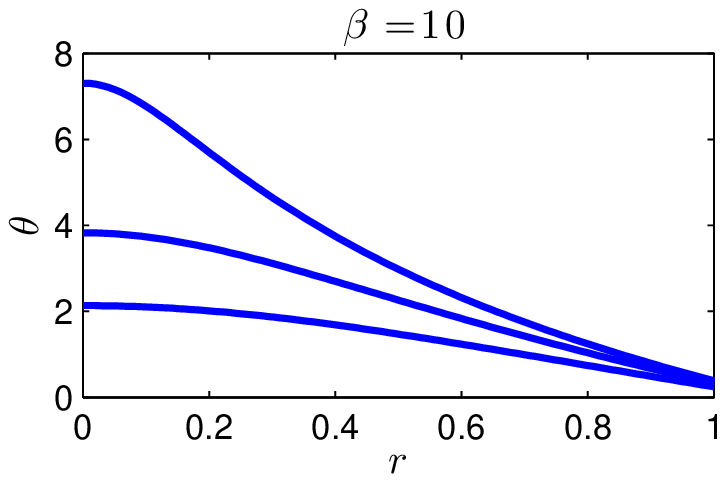}
\caption{\label{flashd_ststV}Steady states of the voltage-controlled problem (\ref{nd2Teq})-(\ref{nd2Tbc1}) for $\beta = 0.1,1,10$, and $\lambda = (0.1,0.5,0.9) \times \lambda_c(\beta)$, as shown by crosses on Figure \ref{flashd_ststV2}.}
\end{figure}

\begin{figure}
\centering
\includegraphics[width=0.32\textwidth]{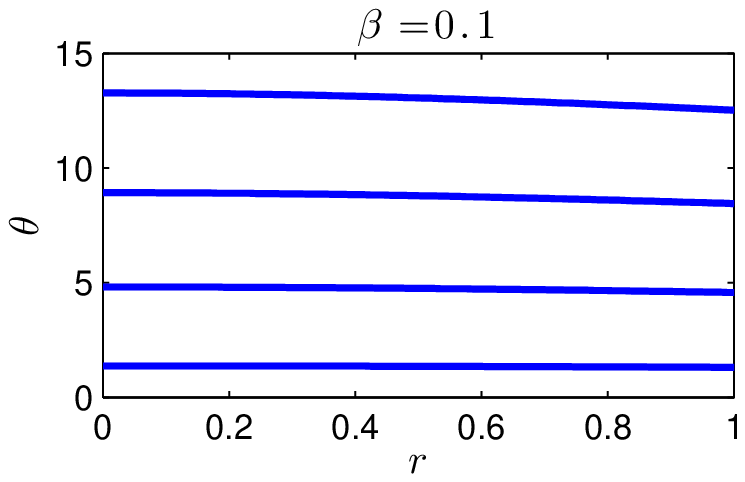}
\includegraphics[width=0.32\textwidth]{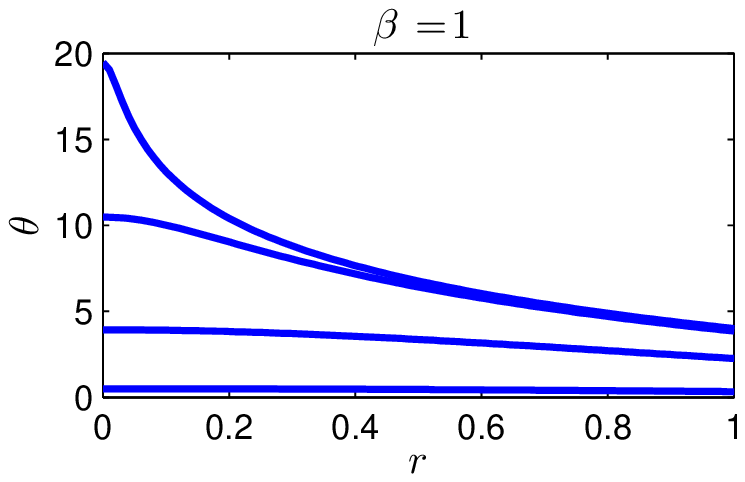}
\includegraphics[width=0.32\textwidth]{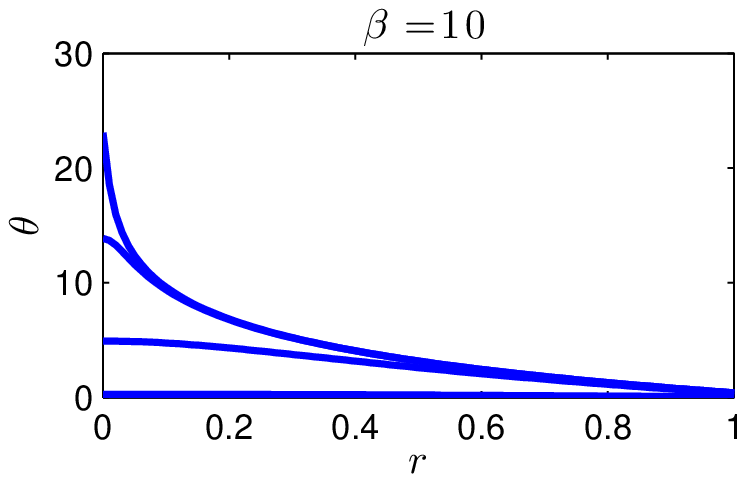}
\caption{\label{flashd_stst}Steady states of the current-controlled problem (\ref{nd2Tbc1})-(\ref{nd2Teq2}) for $\beta = 0.1,1,10$, and $\lambda \mathcal{I}^2 = 1,10^2,10^4,10^6$ (larger $\theta$ for larger $\lambda \mathcal{I}^2$).}
\end{figure}

\begin{figure}
\centering
\includegraphics[width=0.32\textwidth]{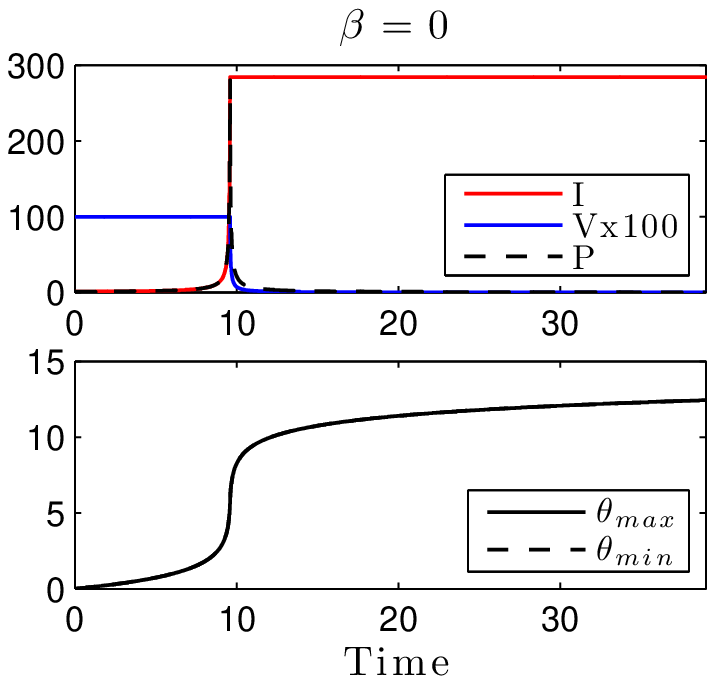}
\includegraphics[width=0.32\textwidth]{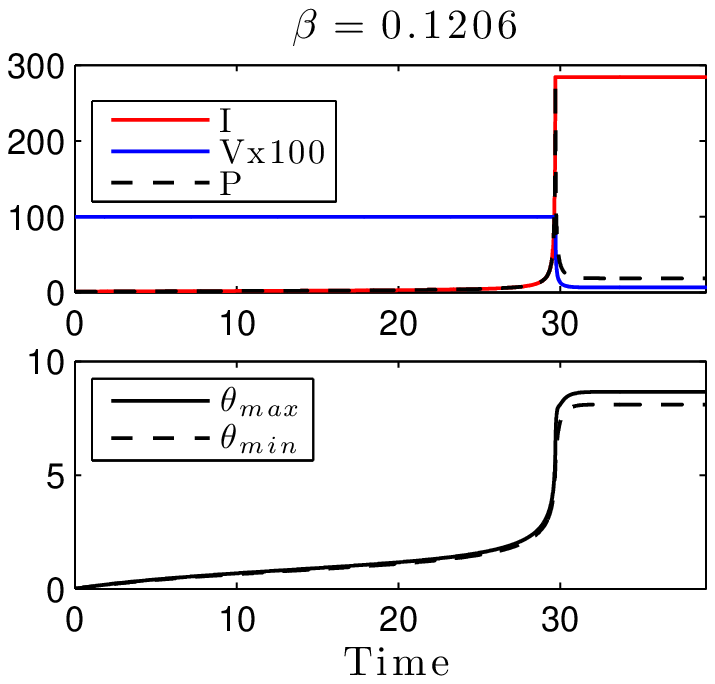}
\includegraphics[width=0.32\textwidth]{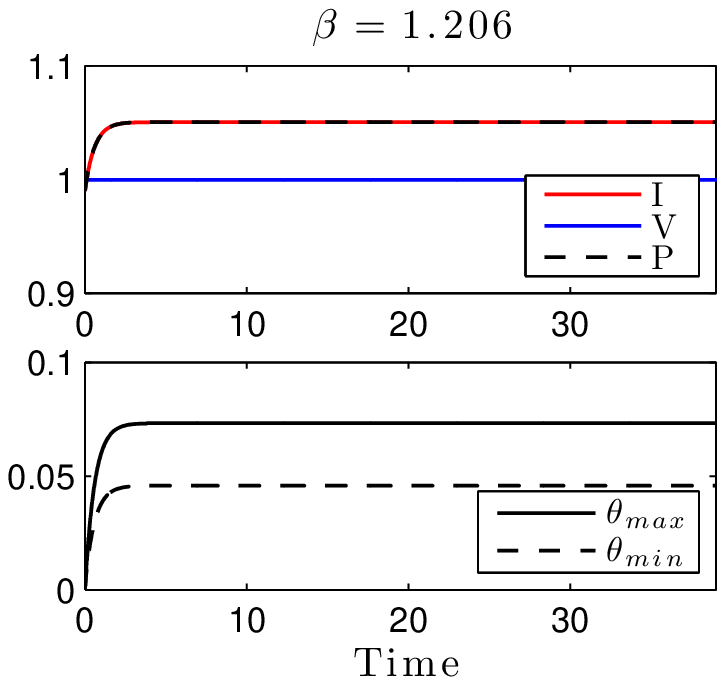}
\caption{\label{flashd}Numerical solutions to (\ref{Veq})-(\ref{nd2Tbc1}) using three different values of $\beta$, and with $\lambda = 0.104$, $\mathcal{I} = 284$.  Top panels show evolution of current, voltage, and power $P = IV$.  Lower panels show the maximum and minimum temperature (at the centre line and outer radius, respectively).  Incipient blow-up, limited by the current limit, is clearly seen in the first two cases.}
\end{figure}

Figure~\ref{flashd_ststV} shows the
steady-state dimensionless temperature profiles
for the voltage-controlled problem
{(\ref{nd2Teq})-(\ref{nd2Tbc1})}
for $\lambda <\lambda_c$, and Figure~\ref{flashd_stst}
shows steady states for the current-controlled problem
{(\ref{nd2Tbc1})-(\ref{nd2Teq2})}. 
An important prediction of Figure~\ref{flashd_stst}
from a practical point of view is that although current control gives
stability there can still be considerable temperature gradients
within the specimen.
Figure~\ref{flashd} shows results from  numerical solutions
to the full problem {(\ref{nd2Teq})-(\ref{nd2Teq2})},
showing the centre ($r=0$) and
surface ($r=1$) temperatures for three different values of $\beta$.
The final case is subcritical, so the switch to current control
is never achieved and the temperature stays close to the furnace temperature.

%%%%%%%%%%%%%%%%%%%%%%%%%%%%%%%%%%%%%%%%%%%%%%%%%%

%%%%%%%%%%%%%%%%%%%%%%%%%%%%%%%%%%%%%%%%%%%%%%%%%%

% % % % % % % % % % % % % % % % % % % % % % % % % % % % % % % % % % % % 

\subsection{Thermally insulating sides}   \label{subsec:electrode}

% % % % % % % % % % % % % % % % % % % % % % % % % % % % % % % % % % % % 

\ijhdel{{ {\bf\{4c\}} This limiting case comes directly from
the experimental parameter values in Table~\ref{parametervalues},
which indicates that the cooling from the sides is very weak,
so that we may set $\beta = 0$. {\it [I'd like to look at
\cite{GSPCJLR} before we refer to it here.]} However, we now allow for the possibility}}
\ijh{Based on the parameter values in Table~\ref{parametervalues}, which indicate that cooling from the sides is rather weak, another limiting case to consider is to set $\beta = 0$.  For more generality, however, we now conversely allow for the possibility}
 of cooling
from the electrodes, with $\alpha >0$.  If the electrodes are
very good conductors of heat they may effectively be held
at the furnace temperature $T_{\infty}$ in which case $\alpha \to \infty$.  

In this case there is nothing to induce any $r$ dependence of the temperature, we have $\theta=\theta(z,t)$, and the electric problem (\ref{ndphieq}) has solution
\be
\phi_z = I \eE^{-\theta}, \qquad 
V = I \int_{-\half}^{\half} \eE^{-\theta}\, {\rm d} z.
\ee
The initial voltage-controlled temperature problem ($V = 1$) is therefore 
\be \label{redeq2}
\theta_t = \delta^2 \theta_{zz} +  \lambda\frac{\eE^{-\theta}}{\left( \int_{-1/2}^{1/2} \eE^{-\theta} \I z \right)^{2}}\, ,
\ee
\be \label{redeqbc}
\theta_z = \mp \alpha \theta \quad {\rm at} \quad z = \pm \half, \qquad
\theta = 0 \quad {\rm at}\quad t = 0.
\ee

Unlike in the previous section, it is now this voltage-controlled regime
which has a non-local term. Moreover, in this case,
the non-local reaction term can lead to blow-up of $\theta$ for all $z$
and hence to blow-up of current ({\cite{Lacey3, Lacey4}};
see below also), and thus there is a switch to current control.
If the end cooling is sufficiently strong, however, the problem may tend
towards a steady state with no blow-up occurring. 

If the switch to current-control is achieved, the equation becomes simply 
\be \label{redeq3}
\theta_t = \delta^2 \theta_{zz} + \lambda\cI^2 \eE^{-\theta} \, , \qquad
\ee
which is a standard reaction-diffusion problem. With the reaction term
decreasing in $\theta$, this problem has a unique equilibrium,
to which all solutions approach for large
time (see \cite{Sattinger}).

\begin{figure}
\centering
\includegraphics[width=0.32\textwidth]{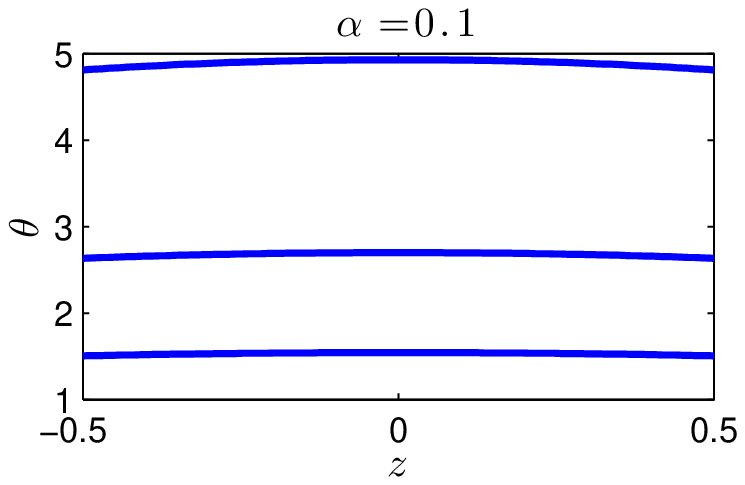}
\includegraphics[width=0.32\textwidth]{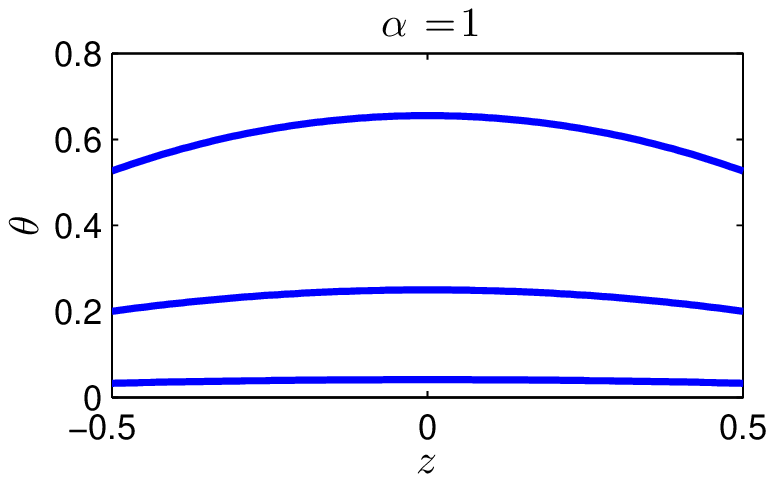}
\includegraphics[width=0.32\textwidth]{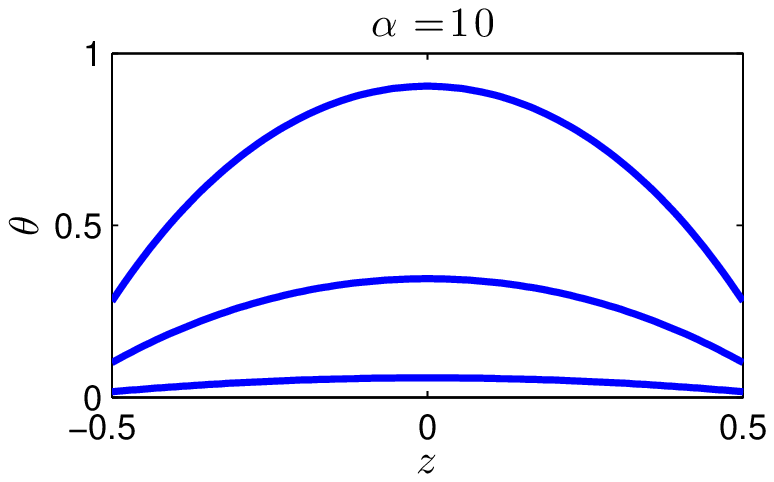}
\caption{\label{flashe_ststV}\ijh{Steady states of the voltage-controlled problem (\ref{redeq2})-(\ref{redeqbc}) for $\alpha = 0.1,1,10$, and $\lambda = (0.1,0.5,0.9) \times \lambda_c(\beta)$, as shown by crosses on Figure \ref{flashd_ststV2}.}}
\end{figure}

For flash to happen in this model therefore requires blow-up to occur
for (\ref{redeq2}), and we determine the conditions under
which this will happen by again seeking a steady-state solution.
If such a solution exists, then blow-up does not occur.
The steady state satisfies
\be \label{redeq4}
\delta^2 \theta_{zz} + \lambda\frac{\eE^{-\theta}}{\left( \int_{-1/2}^{1/2} \eE^{-\theta} \I z \right)^{2}} = 0 \, , \qquad
\theta_z = \mp \alpha \theta \quad {\rm at} \quad z = \pm \half \, ,
\ee
and  the solution is 
\be \label{redsol1}
\theta(z) = 2\ln \left(\frac{\cos az}{\cos \half a} \right) + b
\ee
where the constants $a$ and $b$ are related by
\be
\frac{\lambda}{\delta^2} = 8  \eE^{-b} \sin^2 \half a \, , \qquad
b = \frac{2a}{\alpha} \tan \half a.
\ee
There is a solution for $a$ only if $\lambda  \le \delta^2 \Lambda_c(\alpha)$,
where $\Lambda_c(\alpha)$ is shown in Figure~\ref{flashd_ststV2}.
Solutions of {(\ref{redeq2})-(\ref{redeqbc})} for $\lambda$ less
than this critical value are displayed in Figure~\ref{flashe_ststV}.
It is straightforward to determine that $\Lambda_c \to 8$ for
$\alpha \to \infty$, and $\Lambda_c \sim 2\alpha/\eE$ for $\alpha \to 0$.  

\begin{figure}
\centering
\includegraphics[width=0.32\textwidth]{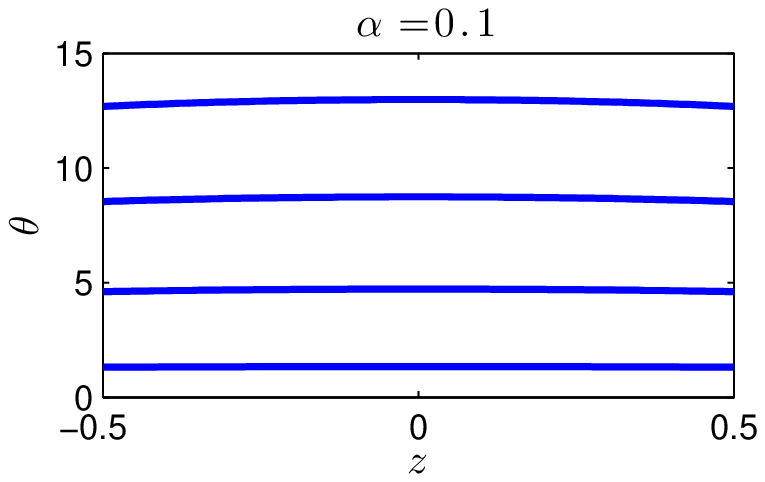}
\includegraphics[width=0.32\textwidth]{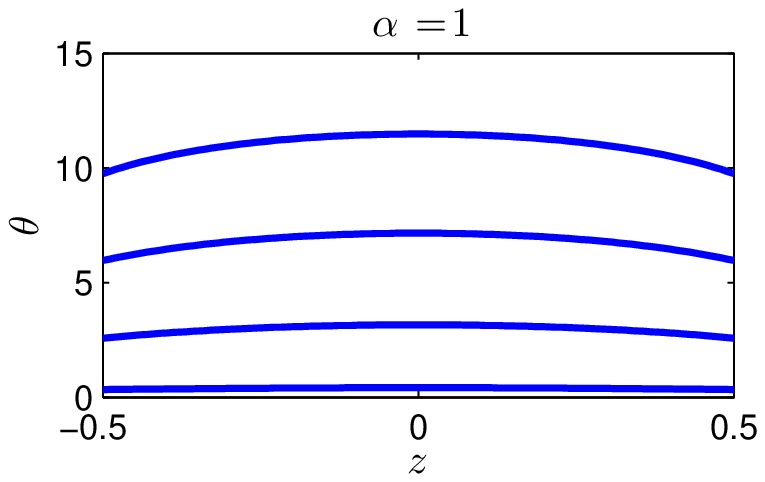}
\includegraphics[width=0.32\textwidth]{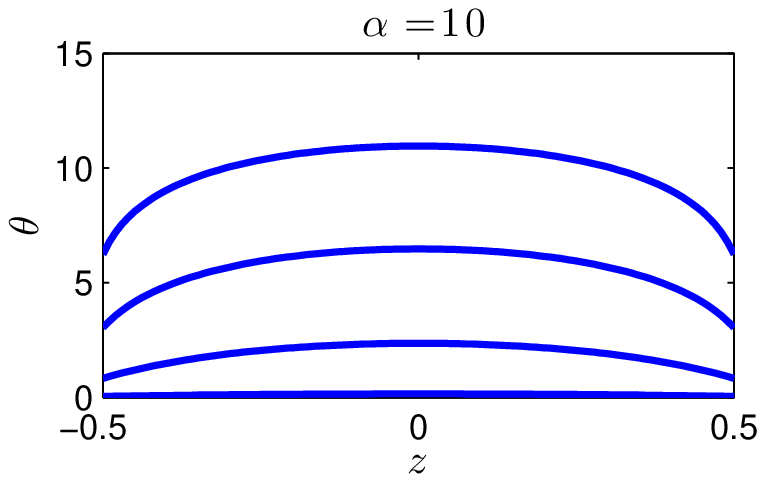}
\caption{\label{flashe_stst}\ijh{Steady states of the current-controlled problem (\ref{redeqbc})-(\ref{redeq3}) for $\alpha = 0.1,1,10$, and $\lambda \mathcal{I}^2 = 1,10^2,10^4,10^6$ (larger $\theta$ for larger $\lambda \mathcal{I}^2$).}}
\end{figure}

Steady states for (\ref{redeq3}) have the same form
as in (\ref{redsol1}), but with
\be
\frac{\lambda \mathcal{I}^2}{\delta^2} = \frac{2a^2}{ \cos^2 \half a} \, e^{b}, \qquad
b = \frac{2a}{\alpha} \tan \half a
\ee
and some sample cases can be seen in Figure~\ref{flashe_stst}.

%%%%%%%%%%%%%%%%%%%%%%%%%%%%%%%%%%%%%%%%%%%%%%%%%%

%%%%%%%%%%%%%%%%%%%%%%%%%%%%%%%%%%%%%%%%%%%%%%%%%%

% % % % % % % % % % % % % % % % % % % % % % % % % % % % % % % % % % % % 

\subsection{Analysis for small \boldmath $\alpha$, $\beta$ and $\lambda$}  \label{subsec:small}

% % % % % % % % % % % % % % % % % % % % % % % % % % % % % % % % % % % % 

If $\alpha = \beta = 0$ in (\ref{ndphieq})-(\ref{ndTbc2}),
the temperature is independent of space, and $\phi = - Vz$,
with $V = \min\left(1,\mathcal{I}/e^{\theta}\right)$.
Concentrating solely on the voltage-controlled regime,
it is clear that the problem exhibits {blow-up} as $t \to 1/\lambda$
according to
\be
\theta = -\log (1-\lambda t).
\ee
Such {blow-up} can only be prevented if the cooling terms are
in fact non-negligible, and $\lambda$ is not too large.
Thus, we consider the case when $\alpha$, $\beta$
and $\lambda$ are all comparably small.  In this case the temperature
is still roughly uniform over the cross-section, and it is appropriate
to integrate over the whole domain, giving
\be \label{fgeq}
\frac{{\rm d}\theta}{{\rm d} t} = f(\theta) - g(\theta),
\ee
where the `bulk' heating and cooling functions are
\be \label{fgdef}
f(\theta) = \lambda \min\left(\eE^{\theta}, \mathcal{I}^2 \eE^{-\theta} \right), \qquad
g(\theta) = 2 \left( \beta + \delta^2 \alpha\right) \theta,
\ee
the cooling terms coming from the side and the ends, respectively.

\begin{figure}
\centering
\includegraphics[width=0.5\textwidth]{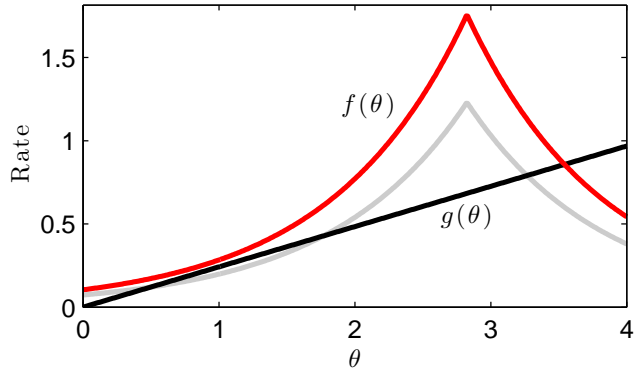}
\caption{\label{flashb_parameters1}\ijh{Heating and cooling functions from (\ref{fgdef}).  The cooling curve $g(\theta)$ is shown in black, and the heating curve $f(\theta)$ is shown in red, for the parameter values in Table~\ref{parametervalues}.  The {steady-state} temperature is given by the intersection of the two curves.  In grey is the heating curve for a smaller (subcritical) value of $\lambda$.} }
\end{figure}

The heat balance (\ref{fgeq}) provides a simple way to understand
the process. Graphs of $f(\theta)$ and  $g(\theta)$ are given in
Figure~\ref{flashb_parameters1} and indicate how the temperature might
evolve towards a steady state
for which $f(\theta) = g(\theta)$.
Flash sintering requires that a sufficiently high temperature is exceeded
before we reach this
steady state. In practice, we expect that this temperature is high enough
that the current-limited regime must always be reached, implying that there
must be no steady state on the voltage-controlled
portion of the heating curve. The absence of such a steady state
is the condition that blow-up would occur, were it not for the current limit.

More precisely, the blow-up condition is that the increasing portion
of the red curve in the Figure~\ref{flashb_parameters1} lies above
the black curve. Otherwise, a steady state will be reached
at too low a temperature.  
This critical condition 
can be expressed as 
\be
\lambda > (2/\eE) (\beta + \delta^2 \alpha) \, , % \Theta > \Theta_{c}(\Lambda),
\label{eq:simplecrit}
\ee 
which is found from the condition that the $f(\theta)$ and $g(\theta)$ curves
meet tangentially.  This is consistent with the limiting behaviour of $\lambda_c(\beta)$ and $\Lambda_c(\alpha)$ \ijh{for small $\beta$ and $\alpha$} found in the previous two sections.

If the current limit $\mathcal{I}$ is too low it is clear
that even the steady state that occurs in the current-limited regime
(which always exists, since $f(\theta) \to 0$ for $\theta \to \infty$)
will be at too low a temperature for sintering.
Thus $\mathcal{I}$ must not be made too small,
but it is otherwise unimportant with regard
to the occurrence of flash.

% % % % % % % % % % % % % % % % % % % % % % % % % % % % % % % % % % % % 

\subsection{High aspect ratio}   \label{subsec:aspect}

% % % % % % % % % % % % % % % % % % % % % % % % % % % % % % % % % % % % 

{We now consider the case when $\delta$ is small, which provides a link between the previous three subsections.  

If $\delta$ is small, but $\alpha$, $\beta$ and $\lambda$ are order one, then the cooling effect of the electrodes is confined to boundary layers at the edges.  In that case the bulk of the sample has temperature independent of $z$ and the analysis of Subsection~\ref{subsec:insulate} applies.  On the other hand, if $\alpha$, $\beta$ and $\lambda$ are also small, the problem is essentially the same as in Subsection~\ref{subsec:small}, and that analysis applies.

If, $\delta$, $\beta$ and $\lambda$ are small, but $\alpha$ is order one (or even infinite, for the case of highly conductive electrodes held at the furnace temperature), another reduction applies.  In that case the temperature is roughly uniform across the radius, and an integral over the cross-section produces the dimensionless model, in the case of voltage control,
\be \label{redeq2:a}
\theta_{\tilde{t}} =
\theta_{zz} +
\Lambda \left( \int_{-1/2}^{1/2} \eE^{-\theta} \I z \right)^{-2} \eE^{-\theta}
- B \theta \, , \qquad
\theta_z = \mp \alpha \theta \quad {\rm at} \quad z = \pm \half,
\ee
where 
\be
\tilde{t} = \delta^2 t , \quad {\Lambda = \lambda/\delta^2
\quad \mbox{ and }} \quad B = 2\beta / \delta^2  \, . 
\ee
After flash (the occurrence of which is discussed below), the current-controlled model
becomes, more simply,
\be \label{redeq3:a}
\theta_{\tilde{t}} =
\theta_{zz} + \Lambda \cI^2 \eE^{-\theta}
- B\theta \, , \qquad
\theta_z = \mp\alpha \theta \quad {\rm at} \quad z = \pm \half,
\ee
which is a standard reaction-diffusion problem. As for (\ref{redeq3}), there is
a unique equilibrium for this current-controlled problem, which all solutions approach for large
time (see \cite{Sattinger}).

To determine whether flash happens, we want to know if blow-up
will occur for the voltage-controlled problem (\ref{redeq2:a}). Once again 
this will occur if $\Lambda$ is sufficiently large, and the critical condition in this case becomes
\be
\Lambda > {\tilde{\Lambda}}_c(\alpha,B) \, . % \Theta > \Theta_c(\Lambda,D).
\ee
To determine this critical value, we seek a steady-state solution, for which
\be \label{redeq4:a}
\theta_{zz} +
\Lambda \left( \int_{-1/2}^{1/2} \eE^{-\theta} \I z \right)^{-2} \eE^{-\theta}
- B\theta = 0 \, , \qquad
\theta_z = \mp \alpha \theta \quad {\rm at} \quad z = \pm \half \, .
\ee
If such a solution exists, then blow-up does not occur. Given the form of (\ref{redeq4:a}),
it is difficult to calculate such a solution analytically and thus determine the critical value
${\tilde{\Lambda}}_c(\alpha,B)$ leading to flash. However,
a lower bound can be obtained by neglecting
either the first or third term in (\ref{redeq4:a}).  

In the case that $B$ is large (as is the case for the data in Table~\ref{parametervalues}), the cooling effect of the electrodes is again confined to diffusive boundary layers at the edges. Over most of the sample the temperature is then approximately the same as in 
Subsection~\ref{subsec:small}, and that earlier analysis shows that $\tilde{\Lambda}_c \approx B/e$.  If $B$ is small, on the other hand, then (\ref{redeq4:a}) reduces to the same problem as in (\ref{redeq2})-(\ref{redeqbc}), and the analysis of Subsection~\ref{subsec:electrode} applies, showing $\tilde{\Lambda}_c \approx \Lambda_c(\alpha)$.}

% \newpage

%%%%%%%%%%%%%%%%%%%%%%%%%%%%%%%%%%%%%%%%%%%%%%%%%%%%%%%%%%%%%%%%

\vspace*{0.5cm}
\setcounter{equation}{0}
\section{Conclusions}              \label{sec:concl}

%%%%%%%%%%%%%%%%%%%%%%%%%%%%%%%%%%%%%%%%%%%%%%%%%%%%%%%%%%%%%%%%

We have described a mathematical model for the process of flash sintering
based on Joule heating, and investigated its behaviour in a simple
cylindrically symmetric geometry.  Numerical and approximate analytical
solutions are consistent with experimental results for flash sintering
\cite{Todd}.  

For certain limiting cases of the parameters, the problem may be reduced
to one spatial dimension and, depending on whether voltage
or current is controlled, the equations can take the form of a local
or non-local reaction diffusion problem.  We have made use
of existing results concerning the solutions of such problems
to infer that flash sintering manifests itself as incipient blow-up behaviour,
and we can identify the controlling processes that determine
whether such behaviour occurs. In all cases this is a competition
between the rate of heating (controlled experimentally by voltage),
and the rate of cooling (predominantly by radiation).
An automatic switch between voltage and current control
in the experiments prevents blow-up from actually occurring,
and the model suggests that once this switch has occurred,
the temperature distribution within the sample evolves towards
a stable steady state. It is clear that the current limit must be
set sufficiently high so that the temperatures needed for rapid sintering
to occur are achieved.  

\begin{figure}
\centering
\includegraphics[width=0.5\textwidth]{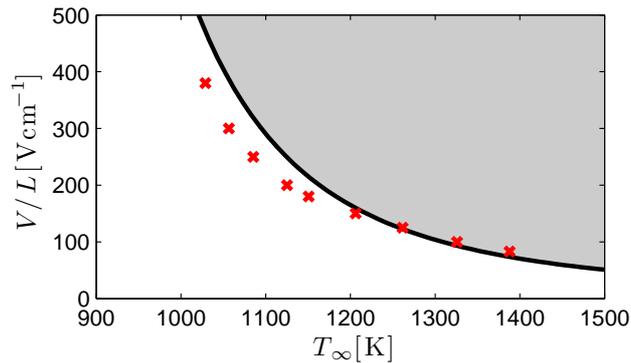}
\caption{\label{flashd_parameters2}Critical condition (\ref{eq:criterion:3}) for {blow-up} of (\ref{nd2Teq})-(\ref{nd2Tbc1}) in terms of dimensional furnace temperature $T_{\infty}$ and electric field $V/L$, all other parameters being held at constant values from Table \ref{parametervalues}.  Flash occurs in the grey region. The red crosses indicate experimental values from \cite{Todd}.}
\end{figure}

The most appropriate of our reduced models to the experimental results
appears to be the case of thermally insulating electrodes,
for which we derived the condition $\lambda > \lambda_c(\beta)$
for a flash to occur, where $\lambda$ is the dimensionless heating parameter,
$\beta$ the dimensionless cooling parameter, and the function
$\lambda_c(\beta)$ is shown in Figure \ref{flashd_ststV2}.
Converted back into dimensional terms, this condition is given as
\begin{equation}
\left( \frac{R^2 AE}{k R_g T_\infty^2} \right) \left( \frac {V}{L} \right)^2
\eE^{-E/(R_g T_\infty)} >  \lambda_c \left( \frac
{R\ijh{h_s} + 4R\epsilon S T_{\infty}^3}k \right).
\label{eq:criterion:3}
\end{equation}
This is expressed as a regime diagram in terms of furnace temperature and voltage
as shown in Figure~\ref{flashd_parameters2}, where the theroretical
predictions are compared with the experimental results of \cite{Todd}.
We note that Figure~\ref{flashd_parameters2}
is similar to Figure 5 of \cite{Todd} but that there
is some systematic overestimate of the critical potential gradient
at lower temperatures. Observe that
condition (\ref{eq:criterion:3}) also demonstrates
the dependence on the geometry
and other parameters.  For instance, increasing the radius of the sample
would make flash more likely.

The model also suggests that before any flash occurs, the temperature
is fairly uniform across the sample, and in that case
Figure~\ref{flashb_parameters1} describes the essential behaviour
quite succinctly. The form of the red heating curve is determined
by the form of
$\sigma(T)$, which was approximated, from an Arrhenius function,
by an exponential.  Literature estimates for the conductivity
and its dependence on temperature vary considerably, however,
so this simple analysis of blow-up may prove useful to investigate
other functional forms for $\sigma(T)$.\footnote{In
finding the critical condition as indicated by Figure~\ref{flashb_parameters1},
we can not only use a full expression for conductivity $\sigma(T)$
but it is also possible to use an un-approximated form of the radiative
cooling law.}

If incipient `blow-up' starts and the current limit is
large enough not to kick in for some time,
{the conclusion that the temperature is almost uniform across the sample
will \ijh{no longer hold true}. This is because the temperature is largest
in the centre of the sample, and leads to faster heating there.
This issue could be investigated further by studying the asymptotics
of near blow-up for small $\beta$, and large $\cI$,
although the local form of such single-point blow-up
has already been looked at,
for example, in \cite{BB, HV1, HV2}: it might be expected,
from consideration of a similarity solution near blow-up at $r=0$,
that the solution of an equation such as (\ref{nd2Teq})
would have a profile at the blow-up time of $\theta \sim - 2 \ln r$ + constant
for $r$ small; however, no suitable similarity solution exists and the usual
profile at blow-up time is somewhat spikier than might be expected,
$\theta \sim - 2 \ln r + \ln (- \ln r)$ + constant.}
% such results indicate a slightly flatter profile than might be expected from a similarity solution. 
From the experimental point of view it is of interest
to know what the maximum temperature is, and what
the temperature differential within the material is.
In this context, it is worth noting that greater temperature variations
across the sample are produced for large $\lambda$, $\mathcal{I}$ and $\beta$
(see the eventual steady states in Figures~\ref{flashd_stst}
and \ref{flashe_stst}). 

Another issue that should be studied further in the future
is the effect of non-uniform cross-section.  By forcing
a non-uniform electric current through the sample,
this may induce a more complex pattern of heating
and enable understanding of experimental problems associated
with partial or incomplete sintering.

%%%%%%%%%%%%%%%%%%%%%%%%%%%%%%%%%%%%%%%%%%%%%%%%%%%%%%%%%%%%%%%%

\vspace*{0.5cm}

\section*{Acknowledgements}
The authors are very grateful to Prof.~J.R.Ockendon~FRS for useful
discussions towards
the preparation of this manuscript.

%%%%%%%%%%%%%%%%%%%%%%%%%%%%%%%%%%%%%%%%%%%%%%%%%%%%%%%%%%%%%%%%

{
}

\end{document}